\documentclass{article} \title{Geometrical construction of quantum groups representations}
\author{Arnal D.\thanks{
Universit\'e de Bourgogne
Laboratoire Gevrey Math\'ematiques Physique
UFR Sciences et Techniques
BP 47870
21078 Dijon Cedex France}~; Bel-Baraka N.$^\ast$\thanks{Universit\'e MohammedV-Agdal Laboratoire de Physique Th\'eorique BP 1014 Rabat Maroc}~; Boukary Baoua O.\thanks{
Universit\'e A.M. de Niamey
Facult\'e des Sciences
D\'epartement de Math\'ematiques
BP~ 10662
Niamey Niger} }
\newtheorem{theo}{Theorem}
\newtheorem{lemme}{Lemma}
\newtheorem{rema}{Remark}

\font\bbfnt=msbm10

\font \tengoth=eufm10
\font \sevengoth=eufm7
\font \fivegoth=eufm5
\newfam\gothfam
\textfont \gothfam=\tengoth
\scriptfont \gothfam=\sevengoth
\scriptscriptfont \gothfam=\fivegoth

\def\BbbN{\hbox{\bbfnt\char'116}}
\def\BbbZ{\hbox{\bbfnt\char'132}}
\def\BbbP{\hbox{\bbfnt\char'120}}

\def\BbbC{\hbox{\bbfnt\char'103}}

\let\frak=\goth

\def\val{\hbox{val}}
\def\vec{\hbox{Vec}}

\def\gb{{}^{\hbox{$G$}}\big/_{\hbox{$B$}}}
\def\gbo{{}^{\hbox{$G_0$}}\big/_{\hbox{$B_0$}}}
\def\gbun{{}^{\hbox{$G_1$}}\big/_{\hbox{$B_1$}}}
\def\gu{{}^{\hbox{$G$}}\big/_{\hbox{$U$}}}
\def\guun{{}^{\hbox{$G_1$}}\big/_{\hbox{$U_1$}}}
\def\guo{{}^{\hbox{$G_0$}}\big/_{\hbox{$U_0$}}}

\begin{document}
\maketitle
\vspace{.8cm}
{\bf Abstract}:
We describe geometrically the classical and quantum inhomogeneous
groups $G_0=(SL(2, \BbbC)\triangleright \BbbC^2)$ and $G_1=(SL(2, \BbbC)\triangleright \BbbC^2)\triangleright \BbbC$ by studying explicitly their shape algebras as a spaces of
polynomial functions with a quadratic relations.%%%%%%%
\section{Introduction}

The problem of describing quantum inhomogeneous group $G$, semi-direct pro-duct
of a semi simple group by an abelian or more generally solvable normal subgroup is still
not totally solved \cite{r5}. The main difficulty is the incompletness of the family
of irreducible finite dimensional representations for these groups.

In the order to overpass this problem, we have either to consider
infinite dimensional irreducible representations or finite dimensional indecomposable representations.\\
This last family of representations is very large and hard to describe (generally
we don't have a classification for such representations). But in the case where the inhomogeneous group is a subgroup of a semi simple one $S$, we can try to restrict ourselves to the family of restrictions to $G$ of irreducible finite dimensional
representations of $S$.

On the other hand, C. Ohn gave a geometrical description of the quantum group
$G = SL(n, \BbbC)$ \cite{r1}, by defining its shape algebra as the space of regular sections
of the line bundles over some sheme in the flag manifold \cite{r4}, or as the space of regular
functions on a manifold $\gu$.

In this paper, we intend to describe explicitly the simplest examples of shape algebra
for an inhomogeneous quantum group by using the geometric approach of C.Ohn.\\
More precisely, we are looking at two subgroups $G_1$ and $G_0$ of $G = SL(3, \BbbC)$,
wich are inhomogeneous of the form $(SL(2, \BbbC)\triangleright \BbbC^2)\triangleright \BbbC$
and $(SL(2, \BbbC)\triangleright \BbbC^2)$.\\
In the Drinfeld-Jimbo quantum universal enveloping algebra  ${\cal U}_q({\frak sl}(3, \BbbC))$, there is a sub-bigebra which is the quantum version ${\cal U}_q(\frak g_1)$ of the classical enveloping algebra of $\frak g_1$. Unfortunately, there is no quantum subalgebra associated to $G_0$ in  ${\cal U}_q({\frak sl}(3, \BbbC))$.\\
We intend to study, geometrically, this phenomena, thus we describe the classical shape algebras
for $G_1$ and $G_0$ as regular functions under the complex manifolds $\guun$ and $\guo$. The first one is dense inside
$\overline{\gu}$ but for the second one $\overline{\guo}$ is a closed submanifold of $\overline{\gu}$ of smaller dimension.\\
In our opinion this is the geometrical presentation of the fact that $G_0$ is not a quantum
subgroup of $G$.\\
We can, finally, quantify the shape algebra of $G_1$ by following the same computation as
for $SL(3)$.
The paper is organized as follows, after recalling our notations and the presentation of
${\cal U}_q(\frak g_1)$ in part 2, we describe explicitly the classical shape algebra for
$G = SL(3, \BbbC)$ as a space of polynomial functions or as a quadratic associative
algebra generated by particular functions $p_i$ and $q_j$ with explicit quadratic relations.\\
In part 4, we recall the C. Ohn construction for the quantum shape algebra of $G$. As a result, this algebra is still quadratic, generated by the $p_i$ and $q_j$ with explicit deformed
relations.\\
In the two last part, we describe first the classical shape algebras for $G_0$ and $G_1$
as algebras of polynomial functions, then the quantum shape algebra for $G_1$, as an explicit
quadratic associative algebra generated by $p_i$, $q_j$ and $1/{q_3}$.

\section{Notations and preliminaries}

In this paper, we shall consider the following Lie groups:
$$
G=SL(3,\BbbC)=\{
\begin{array}{l}
g =\left(\begin{array}{ccc}
x_1&y_1&z_1\\
x_2&y_2&z_2\\
x_3&y_3&z_3\\
\end{array}\right)
\end{array}, \quad \  \det g = 1 \}
$$
$$
G_1=\left( SL(2,\BbbC) \triangleright \BbbC^2\right) \triangleright \BbbC =\{
\begin{array}{l}
g_1=\left(\begin{array}{ccc}
x_1&y_1&0\\
x_2&y_2&0\\
x_3&y_3&z_3\\
\end{array}\right)
\end{array}, \quad \  \det g_1 = 1\}
$$
and
$$
G_0=SL(2,\BbbC)\triangleright\BbbC^2 = \{
\begin{array}{l}
g_0=\left(\begin{array}{ccc}
x_1&y_1&0\\
x_2&y_2&0\\
x_3&y_3&1\\
\end{array}\right)
\end{array}, \quad \  \det g_0 = 1\}
.$$
$G$ being simple, there exists quantum versions of $G$. Especially here, we consider the well
 known Drinfeld-Jimbo quantum universal enveloping algebra ${\cal U}_q({\frak g})$ for
 ${\frak g}= {\frak sl} (3,\BbbC)$, it is defined by its generators
$K^{\pm}_1, K^{\pm}_2, X_1, X_2, Y_1, Y_2$ and the relations:
$$
\begin{array}{ccc}
K_1K_1^{-1}=K_1^{-1}K_1=1\\
\vspace{.2cm}
K_2K_2^{-1}=K_2^{-1}K_2=1 &
\vspace{.2cm}
K_1K_2=K_2K_1& (*)\\
\vspace{.2cm}
K_1 X_1 K_1^{-1} = q^2 X_1\hfill &
K_1 Y_1 K_1^{-1} = q^{-2} Y_1\hfill&\hskip5mm
K_1 X_2 K_1^{-1} =q^{-1}X_2\\
\vspace{.2cm}
K_1 Y_2 K_1^{-1} =q Y_2\hfill &
K_2 X_1 K_2^{-1} = q^{-1} X_1\hfill&
 K_2 Y_1 K_2^{-1} = q  Y_1\\
\vspace{.2cm}
K_2 X_2 K_2^{-1} =q^2 X_2\hfill&
K_2 Y_2 K_2^{-1} =q^{-2} Y_2\hfill&
X_2Y_1-Y_1X_2 =0\\
\vspace{.2cm}
X_1Y_2-Y_2X_1=0\hfill&
X_1Y_1\!-\!Y_1X_1 \!=\!\frac{K_1-K_1^{-1}}{q-q^{-1}}&\hskip4mm
X_2Y_2\!-\!Y_2X_2\! =\!\frac{K_2-K_2^{-1}}{q-q^{-1}}\\
\end{array}
$$
and
$$
\begin{array}{cc}
X_1^{2}\! X_2\!-\!(q\!+\!q^{-1})\!X_1\!X_2\!X_1\!+\!X_2\!X_1^2\!=\!0
\hskip3mm&\hskip3mm Y_1^{2}\! Y_2\!-\!(q\!
+\!q^{-1})\!Y_1\!Y_2\!Y_1\!+\!Y_2\!Y_1^2\!=\!0\\
\vspace{.3cm}
X_2^{2}\! X_1\!-\!(q\!+\!q^{-1})\!X_2\!X_1\!X_2\!+\!X_1\!X_2^2\!=\!0
\hskip3mm&\hskip3mm
Y_2^{2}\! Y_1\!-\!(q\!+\!q^{-1})\!Y_2\!Y_1\!Y_2\!+\!Y_1\!Y_2^2\!=\!0\\
\end{array}
$$
The coalgebra structure on ${\cal U}_q({\frak g})$
 is defined by the following coproduct:
$$
\begin{array}{ccc}
\Delta\! K^{\pm 1}_1\!=\!K_1^{\pm 1}\!\otimes\! K_1^{\pm 1}
\hskip2mm&\hskip2mm
\Delta\! X_1=X_1\!\otimes\! 1\!+ \!K_1\!\otimes \!X_1
 \hskip2mm&\hskip2mm
\Delta\! Y_1\!=\!Y_1\!\otimes\! K_1^{-1}\!+\! 1\!\otimes\! Y_1\\
\vspace{.3cm}
\Delta\! K^{\pm 1}_2\!=\!K_2^{\pm 1}\!\otimes\! K_2^{\pm 1}
\hskip2mm&\hskip2mm
\Delta\! X_2\!=\!X_2\!\otimes\! 1\!+\! K_2\!\otimes\! X_2&
\Delta\! Y_2\!=\!Y_2\!\otimes\! K_2^{-1}\!+\! 1\!\otimes\! Y_2\\
\end{array}
$$
The generators $K^{\pm 1}_i$ are formally identified with $e^{\pm tH_i}$
 if $H_1$ and $H_2$ are the usual basis for the Cartan subalgebra :
$$
\begin{array}{l}
H_1= \left(\begin{array}{ccc}
1&0&0\\
0&-1&0\\
0&0&0\\
\end{array}\right)
\hskip1cm
H_2= \left(\begin{array}{ccc}
0&0&0\\
0&1&0\\
0&0&-1\\
\end{array}\right)
\end{array}
$$
$X_i$ and $Y_i$ are the root vectors associated to simple roots :
$$
\begin{array}{l}
X_1=\left(\begin{array}{ccc}
0&1&0\\
0&0&0\\
0&0&0\\
\end{array}\right)
\hskip0,5cm
X_2=\left(\begin{array}{ccc}
0&0&0\\
0&0&1\\
0&0&0\\
\end{array}\right)
\end{array}
$$
$$
\begin{array}{l}
Y_1=\left(\begin{array}{ccc}
0&0&0\\
1&0&0\\
0&0&0\\
\end{array}\right)
\hskip0,5cm
Y_2=\left(\begin{array}{ccc}
0&0&0\\
0&0&0\\
0&1&0\\
\end{array}\right)
\end{array}
$$
and $t$ and $q$ are formally related by $e^t=q$.\par
 Since the Lie algebras $\frak g_0$ and $\frak g_1$ of $G_0$
and $G_1$ are subalgebras of $\frak g$, the classical universal
 enveloping algebras
${\cal U}({\frak g}_0)$ and ${\cal U}({\frak g}_1)$ are subalgebras of ${\cal U}({\frak g})$.
 Unfortunately, this does not hold at the quantum level.
\begin{lemme}

\

\noindent
${\cal U}_q({\frak g}_1)$ is a sub-bigebra of ${\cal U}_q({\frak g})$
but $\Delta ({\cal U}_q({\frak g}_0))$ is not included in
 ${\cal U}_q({\frak g}_0)\otimes{\cal U}_q({\frak g}_0)$.
\end{lemme}
Indeed, the bigebra ${\cal U}_q({\frak g}_1)$ can be written as :
$$
{\cal U}_q({\frak g}_1)=\ {}^
{\displaystyle T(K_1^{\pm 1},K_2^{\pm 1},X_1,Y_1,Y_2)}
\bigg\slash_{\displaystyle T(K_1^{\pm 1},K_2^{\pm 1},X_1,Y_1,Y_2)
\cap I}
$$
where $T(K_1^{\pm 1},K_2^{\pm 1},X_1,Y_1,Y_2)$
is the tensor algebra generated by
$K_1^{\pm 1},K_2^{\pm 1},X_1,$ $Y_1$ and $Y_2$ and $I$ is the ideal given by
 relations $(*)$ in which the generator $X_2$ does not arise.
Moreover, the coproduct $\Delta$ of ${\cal U}_q({\frak sl}(3,\BbbC))$
 can be restricted to
${\cal U}_q({\frak g}_1)$ as :
$$
{\Delta_\vert}_{{\cal U}_q({\frak g}_1)}:
{\cal U}_q({\frak g}_1)\rightarrow {\cal U}_q({\frak g}_1)\otimes {\cal U}_q({\frak g}_1).
$$
Neverthless, the set of generators of ${\cal U}_q({\frak g}_0)$
 does not contain $K_2$. Since $K_2$ arise in the expression of $\Delta Y_2$,
 then ${\cal U}_q({\frak g}_0)$ is not included in
${\cal U}_q({\frak g}_0)\otimes {\cal U}_q({\frak g}_0)$.\vskip2mm

Due to this fact, many authors, attempting to describe inhomogeneous
 quantum goups like
$SL(2,\BbbC)\triangleright\BbbC^2$, prefer to add a central dilatation
 (the $K_2$ element here) to build their model (see \cite{r5} for instance). \par
In this paper, we want to give another point of view more geometrical
for this phenomena. Our starting point is the description, by C. Ohn \cite{r1}, of quantum
$SL(3,\BbbC)$ and its shape algebra.

\section{Borel-Weil-Bott theorem for $SL(3,\BbbC)$}
\label{sec}
To describe $SL_q(3,\BbbC)$, we need an explicit realization for each
irreducible unitary representation of $SL(3,\BbbC)$ and their tensor product.
 The unitary irreducible representation of $SL(3,\BbbC)$ are acting on the
 space of sections of line bundles over $\gb$, where $B$
is the Borel subgroup \cite{r4}:$$
B=\{
\begin{array}{l}
b =\left(\begin{array}{ccc}
a_1&b_1&c_1\\
0&b_2&c_2\\
0&0&c_3\\
\end{array}\right)
\end{array}, \quad \  \det b = 1 \}
$$\\
The characters $\chi_{n_1,n_2}$ of $B$ have the following form :
$$
\chi_{n_1,n_2}(b) =a_1^{-n_1}c_3^{n_2}=a_1^{-(n_1+n_2)}b_2^{-n_2}
\hskip1cm (n_1,n_2\in\BbbZ).
$$
We denote the corresponding principal bundle by:
$$
G\times_{\chi_{n_1,n_2}}\BbbC\rightarrow\gb
$$
(an element of  $G\times_{\chi_{n_1,n_2}}\BbbC$ is an equivalence class
$[g,z]=[gb,\chi_{n_1,n_2}(b)^{-1}z]$ in $G\times \BbbC$).\\
$\gb$ is the flag manifold $D$ of $\BbbC^3$:
$$
D=
\begin{array}{l}\{
\BbbC\!\left[\begin{array}{c}x_1\\x_2\\x_3\end{array}\right],
\BbbC\!\left[\begin{array}{c}x_1\\x_2\\x_3\end{array}\right]+
\BbbC\!\left[\begin{array}{c}y_1\\y_2\\y_3\end{array}\right]
\}
\simeq \{ [p],[q],\hskip0.1cm p.q= p_1q_1+p_2q_2+p_3q_3=0 \}
\end{array}
$$
where the elements $[p]$ and $[q]$ of $\BbbP (\BbbC^3)$ are respectively
the line through\par
$p~=~\left[\begin{array}{c}p_1=x_1\\p_2=x_2\\p_3=x_3\end{array}\right]$
and the line trough
$q=\left[\begin{array}{c}q_1=x_2y_3-x_3y_2\\q_2=x_3y_1-x_1y_3\\
q_3=x_1y_2-x_2y_1\end{array}\right]$ in $\BbbC^3$.
The space of holomorphic sections of the line bundle
 $G\times_{\chi_{n_1,n_2}}\BbbC\rightarrow\gb$ is non trivial if and only if
 $n_1\geq 0$ and $n_2\geq 0 $.
A section is a regular homogeneous function $f$ from $G$ to $\BbbC$ such that :
$$
f(gb)=\chi_{n_1,n_2}(b)^{-1}f(g),\hskip0,5cm \forall g\in G, \forall b\in B.
$$
Let $U$ be the unipotent subgroup :
$$
U=\{u=\left(\begin{array}{ccc}
1&b_1&c_1\\
0&1&c_2\\
0&0&1\\
\end{array}\right)\}.
$$
We can describe the homogeneous space $\gu$ as the submanifold of $\BbbC^6$
 defined by :
$$
\gu\simeq\{ p\in \BbbC^3\backslash\{0\},~~q\in \BbbC^3\backslash\{0\},~~
 p.q= 0\}.
$$
$\gu$ is thus an affine submanifold of $\BbbC^6$, we shall consider its closure
 $\overline{\gu}$ in the usual embeding of $\BbbC^6$ into $\BbbP(\BbbC^7)$:
$$
\overline{\gu }
=\{\left[\begin{array}{l}p\\q\\1\end{array}\right]\in\BbbP(\BbbC^7),\ p.q=0\}
$$
\begin{lemme}[Extension of sections]

\

\noindent
For each $n_1, n_2$ in $\BbbN$, the sections of
$G\times_{\chi_{n_1,n_2}}\BbbC\rightarrow\gb$ can be viewed as regular
functions on $\overline{\gu}$.
\end{lemme}\vskip2mm

\noindent
{\bf Proof:}
Since $\chi_{n_1,n_2}$ is trivial on $U$, $f$ gives rise to a function, still
 denoted $f$ on $\gu$.\\
 If $p_1\not= 0$ and $q_3\not=0$, then
$$
\left(\begin{array}{ccc}
x_1&y_1&z_1\\
x_2&y_2&z_2\\
x_3&y_3&z_3\\
\end{array}\right)
=\begin{array}{l}
\left(\begin{array}{ccc}
p_1&0&0\\
p_2&\frac{q_3}{p_1}&0\\
p_3&\frac{-q_2}{p_1}&\frac{1}{q_3}\\
\end{array}\right).
\left(\begin{array}{ccc}
1&\frac{y_1}{x_1}&\frac{z_1}{x_1}\\
0&1&\frac{z_2x_1-x_2z_1}{q_3}\\
0&0&1\\
\end{array}\right)
\end{array}
.$$
Let $f$ be a section of our line bundle, then, after division by $\det^k (g)$ with a well choice of $k$,$f$ is polynomial homogeneous in
 the $(x,y,z)$ variables with degree $n_1+n_2$ in $x$, $n_2$ in $y$, $0$
 in $z$ or :
$$
f(x,y,z)=\frac{P_1(p,q)}{p_1^{n_2}}\hskip0,5cm (p_1\not = 0),
$$
$P_1$ being polynomial in $p$ and $q$, homogeneous with degree $n_1+n_2$ in $p$
  and homogeneous  with degree $n_2$ in $q$.\\
Similarly, on the open subset $p_2\not=0$, $p_3\not =0$, we can write :
$$
f(x,y,z)=\frac{P_2(p,q)}{p_2^{n_2}}\hbox{ and }~~~f(x,y,z)=\frac{P_3(p,q)}{p_3^{n_2}}
$$
On the other hand, the ideal ~~$(pq)$~~ generated by ~~$pq=p_1q_1+p_2q_2+p_3q_3$~~
 is prime.
Indeed, if $P$ and $Q$ are polynomials such that:
$$
PQ=(pq).R,
$$
then, by division in $\BbbC(q_1)[p_1,p_2,p_3,q_2,q_3]$, we get polynomial
functions $T, S, T', S'$ such that:
$$
q_1^lP=(pq)T+S \mbox~~{\rm and }~~~q_1^{l'}Q=(pq)T'+S'
$$
with $\deg_{p_1}S=\deg_{p_1}S'=0$. Or :
$$
q_1^{l''}(pq)R=(pq)((pq)TT'+TS'+T'S)+ SS'
$$
Thus $SS'=(pq)R'$, this implies $R'=0$ and $S$ or $S'=0$.
If for instance $S=0$ then
$$
q_1^lP=(pq)T
$$
Then $\val_{q_1}(T)\geq l$ or $T=q_1^lT'$, $P=(pq)T'$ is in $(pq)$.\\
Now by the Nullstellensatz, our equation:
$$f=\frac{P_1}{p_1^{n_2}}=
\frac{P_2}{p_2^{n_2}}~~~~p_1\not=0, ~~~~p_2\not=0 ~~~~ p.q=0$$
 can be written as:
$$
p_2^{n_2}P_1 -p_1^{n_2}P_2=(pq)Q.
$$
Thus $\val_{(p_1,p_2)}Q\geq n_2$ or $Q=\sum_{j=0}^{n_2}p_1^{n_2-j}p_2^jQ_j$.
From that we get :
$$
p_2^{n_2}(P_1-(pq)Q_{n_2})-p_1^{n_2}(P_2-(pq)Q_0)= (pq)p_1p_2
\sum_{j=1}^{n_2-1}p_1^{n_2-j-1}p_2^{j-1}Q_j.
$$
Thus
$$
\val_{p_1}(P_1-(pq)Q_{n_2})\geq 1,\hskip0,5cm \val_{p_2}(P_2-(pq)Q_0)\geq 1
$$
And
$$
P_1-(pq)Q_{n_2}=p_1P_1', \hskip0,5cm P_2-(pq)Q_0=p_2P_2'
$$
then
$$f=\frac{P'_1}{p_1^{n_2-1}}=\frac{P'_2}{p_2^{n_2-1}}$$ and  by induction
there exists a polynomial function $P$ in the variables $p,q$ such that
 $f=P$ on $\gu$.

We shall call ${\cal O}(\overline{\gu})$ the shape algebra of the classical
 group $SL(3,\BbbC)$. The multiplication in this algebra:
$$
{\cal O}(\overline{\gu})\otimes {\cal O}(\overline{\gu})\rightarrow
{\cal O}(\overline{\gu})
$$
is the dual form of the classical comultiplication
 $\Delta$ on ${\cal U}({\frak sl}(3))$.
If $u$ belongs to ${\cal U}({\frak sl}(3))$ and if $f_1,f_2$ belongs to
${\cal O}(\overline{\gu})$, we can put  :
$$
\Delta u (f_1 \otimes f_2)(e,e) =u(f_1f_2)(e)
$$
(in fact $\Delta X =1\otimes X+X\otimes 1$ for any $X$ in
 ${\frak sl}(3,\BbbC))$ and $\Delta$ is a morphism from ${\cal U}({\frak sl}(3))$ into
${\cal U}({\frak sl}(3))\otimes {\cal U}({\frak sl}(3))$.

We can also see, as C. Ohn did in \cite{r1}, the space
$H^0(\gb,\chi_{n_1,n_2})$ of polynomial functions homogeneous with degree $n_1$
 in $p$ and $n_2$ in $q$ as the dual
$V_{n_1\varpi_1+n_2\varpi_2}$ of the ``algebraic'' space
$
V^{n_1\varpi_1+n_2\varpi_2}$. If $\varpi_1$ and $\varpi_2$ stands for the fundamental weights,
$V^{n_1\varpi_1+n_2\varpi_2}$ is the carrying space of irreducible representation of
 $SL(3,\BbbC)$ with highest weight $\lambda=n_1\varpi_1+n_2\varpi_2$. If $V^1$ is $\BbbC^3$ with canonical basis $e_1,e_2,e_3$, then $V^{n_1\varpi_1+n_2\varpi_2}$ is explicitly realized as the submodule of $(V^1)^{\otimes n_1}\otimes (V^1\wedge V^1)^{\otimes n_2}$ generated by the highest
weight vector $v_\lambda =(e_1)^{\otimes n_1}\otimes (e_1\wedge e_2)^{\otimes n_2}$.
Denote $V^2$ the space $V^1\wedge V^1$, we get an inclusion mapping:

$$ V^{n_1\varpi_1+n_2\varpi_2}\subset(V^1)^{\otimes n_1}\otimes (V^2)^{\otimes n_2}.
$$
The natural identification between $H^0(\gb,\chi_{n_1,n_2})$ and
$V_{n_1\varpi_1+n_2\varpi_2}$ being
$$
\varphi\in V_{n_1\varpi_1+n_2\varpi_2}\mapsto f
$$
such that $f(g)=\varphi(g.v_\lambda)$,~for any $g$ in  $G$.\\
 The multiplication $m$
 is thus the transposition of the family of injections :
$$
V^{\lambda_1+\lambda_2}\rightarrow V^{\lambda_1}\otimes V^{\lambda_2}.
$$
In fact as an algebra, ${\cal O}(\overline{\gu})$ is generated by $V_1$
 and $V_2$ the duals of the fundamental representation of $ SL(3,\BbbC)$
 i.e. by the linear functions $p_1,p_2,p_3,q_1,q_2,q_3$ and sixteen
 quadratic relations :
$$
\begin{array}{ccc}
p_ip_j=p_jp_i\hskip2mm (i<j)\hskip5mm\hfill& q_jq_i=q_iq_j \hskip2mm (i<j)&\\
p_iq_j=q_jp_i \hskip2mm (i\not =j)\hfill
& p_iq_i=q_ip_i \hfill&\hfill\hskip5mm p.q=0\\
\end{array}
$$
We can replace the four last relations by :
$$
\begin{array}{cc}
p_1q_1=q_1p_1\hfill& \hskip5mm p_3q_3=q_3p_3\hfill\\
p_2q_2+q_1p_1+q_3p_3=0\hskip5mm&\hskip5mm q_2p_2+p_1q_1+p_3q_3=0
\end{array}
$$

\section{The C. Ohn construction}
C. Ohn gave a more geometrical construction for the shape algebra\cite{r1}.
 Let us recall quickly here, in the $SL(3,\BbbC)$ case, his construction.\\
 The shape algebra is generated by $V_1\oplus V_2$
 (the linear polynomial functions in $p_1,p_2,p_3,q_1,q_2,q_3$).
 If $i$ and $j$ are in $\{1,2\}$,
 we consider the tensor products $V^i\otimes V^j$.

Let $ V^{ij}$ be the irreducible module with highest weight $\varpi_i+\varpi_j$.
 We define an explicit injective map :
$$
V^{ij}\longrightarrow V^i\otimes V^j
$$
by constructing a system of vectors  $e^{ij}_C$ generating $V^{ij}$.
We consider all the orthocell $C$ for ${\frak sl}(3)$:
$C$ is the right coset in the weyl group $W$ of ${\frak sl}(3)$, for a subgroup
 $\Gamma$ generated by a set $A$ of pairwise commuting reflexions, ($A=\emptyset$
 or $A=\{ s_{\alpha_i}\}\ (i=1,2)$).
Among all these orthocells, we select the small and $ij$-effective ones,
(See \cite{r1} for explicit definition).\\ In the  case of $SL(3)$ we get
 fourteen small orthocells. Identifying the Weyl group as $S_3$,
 the six `` trivials'' are:
$$
\begin{array}{ccc}
C_1^0=\left\{\left[ 123\right]\right\}\hskip1cm;\hskip1cm
&C_2^0=\left\{\left[ 132\right]\right\}\hskip1cm;\hskip1cm
&C_3^0=\left\{\left[ 213\right]\right\}\vspace{0.3cm}\\
C_4^0=\left\{\left[ 231\right]\right\}\hskip1cm;\hskip1cm
&C_5^0=\left\{\left[ 312\right]\right\}\hskip1cm;\hskip1cm
&C_6^0=\left\{\left[ 321\right]\right\}.\\
\end{array}
$$
and the eight ``non trivials'' are:
$$
\begin{array}{cc}
&\hskip-1cm C_1=\left\{\left[ 123\right],\left[ 213\right]\right\};\hskip.3cm
C_2=\left\{\left[ 132\right],\left[ 231\right]\right\};\hskip.3cm
C_3=\left\{\left[ 312\right],\left[ 321\right]\right\}\vspace{0.2cm}\\

&\hskip-1cm C_4=\left\{\left[ 123\right],\left[ 132\right]\right\};\hskip.3cm
C_5=\left\{\left[ 213\right],\left[ 312\right]\right\};\hskip.3cm
C_6=\left\{\left[ 231\right],\left[ 321\right]\right\}\vspace{0.2cm}\\

&\hskip-1cm C_7=\left\{\left[ 132\right],\left[ 312\right]\right\};\hskip.3cm
C_8=\left\{\left[ 213\right],\left[ 231\right]\right\}.\vspace{0.2cm}\hfill\\
\end{array}
$$
For any $i,j$ the six trivial orthocells are $ij$-effective.
 Amid the others, we keep only $C_1,C_2,C_5,C_6,C_7$ for $i=j=1$,
   $C_2,C_3, C_4,C_5,C_8$ For $i=j=2$ and   $C_2$, $C_3$ For $i=1$ and $j=2$
 (or $i=2$ and $j=1$). To each orthocell $C$,
 we associate a vector $e_C^{ij}$ of $V^i\otimes V^j$ by the following rule :\\
First we realize the Weyl group $W$ as permutations matrices in $SL(3)$, then
 we put $w.e^{(i)}=e^i_w$ if $e^{(i)}$ is the highest weight vector for $V^i$
 and finally :
$$
e^{ij}_C =\sum_{L\subset A} e^i_{s_{\bar L}w}\otimes e^j_{s_{ L}w}
$$
where $C=\{ s_{\alpha_i}w , \alpha_i\in A\}$,~ $\bar L =A\backslash L$ and $s_L$
 ~is the product of $s_{\alpha_i}$ for $\alpha_i$ in $L$. In the $SL(3)$
 case we get the following vectors :
$$
\begin{array}{cc}
 e^{11}_{C^0_1}=e^{11}_{C^0_2}=e_1\otimes e_1\hfill
&\hskip3mm  e^{11}_{C^0_3}=e^{11}_{C^0_4} =e_2\otimes e_2
\hfill\vspace{.2cm}\\
 e^{11}_{C^0_5} =e^{11}_{C^0_6}=e_3\otimes e_3\hfill
&\hskip3mm e^{11}_{C_1}=e^{11}_{C_2}= e_1\otimes e_2 +q e_2\otimes e_1
\vspace{.2cm}\\
 e^{11}_{C_5}=e^{11}_{C_6} = e_2\otimes e_3+ q e_3\otimes e_2\hfill
&\hskip3mm  e^{11}_{C_8} =e_1\otimes e_3 + q e_3\otimes e_1
\hfill\vspace{.6cm}\\
 e^{22}_{C^0_1}=e^{22}_{C^0_3}=(e_1\wedge e_2)\otimes (e_1\wedge e_2)\hfill
&\hskip3mm  e^{22}_{C^0_5}=e^{22}_{C^0_2}
=(e_1\wedge e_3)\otimes (e_1\wedge e_3)\hfill \vspace{.2cm}\\
e^{22}_{C^0_4}=e^{22}_{C^0_6}=(e_2\wedge e_3)\otimes (e_2\wedge e_3)
&\hfill \\
\end{array}
$$
$$
\begin{array}{c}
e^{22}_{C_2}=e^{22}_{C_3}=(e_1\wedge e_3)\otimes (e_2\wedge e_3)+ q
(e_2\wedge e_3)\otimes (e_1\wedge e_3)\hfill \vspace{.2cm}\\
e^{22}_{C_4}=e^{22}_{C_5}=(e_1\wedge e_2)\otimes (e_1\wedge e_3)+ q
(e_1\wedge e_3)\otimes (e_1\wedge e_2)\hfill \vspace{.2cm}\\
e^{22}_{C_8}=(e_2\wedge e_1)\otimes (e_2\wedge e_3)+ q
(e_2\wedge e_3)\otimes (e_2\wedge e_1)\hfill
\end{array}
$$
\vspace{.6cm}
$$
\begin{array}{c}
e_{C_1^0}^{12}= e_1\otimes (e_1\wedge e_2)
\hskip.5cm \hskip.5cm
e_{C_4^0}^{12}= e_2\otimes (e_2\wedge e_3)\vspace{.2cm}\\
e_{C_2^0}^{12}= e_1\otimes (e_1\wedge e_3)
\hskip.5cm \hskip.5cm
e_{C_5^0}^{12}= e_3\otimes (e_1\wedge e_3)\vspace{.2cm}\\
e_{C_3^0}^{12}= e_2\otimes (e_1\wedge e_2)
\hskip.5cm \hskip.5cm
e_{C_6^0}^{12}= e_3\otimes (e_2\wedge e_3)\vspace{.2cm}\\
e_{C_2}^{12}= e_1\otimes (e_2\wedge e_3) +q e_2\otimes (e_1\wedge e_3)
\hskip.2cm \hskip.2cm
e_{C_5}^{12}= e_2\otimes (e_3\wedge e_1) +q e_3\otimes (e_2\wedge e_1)\\
\end{array}.
$$
With these notations $V^{ij}$ is linearly generated by the vectors
 $e^{ij}_C$ where $C$ is small and $ij$-effective \cite{r1}.
\begin{rema}
This construction can be compared with the Demazure's one of a basis
  for an irreducible module for a simple Lie algebra \cite{r6},in fact, the explicit C. Ohn
construction gives a generating system only if the highest weight has the form
$\lambda_1+\lambda_2$ where $\lambda_1$ and $\lambda_2$ are fundamental.
\end{rema}
\vspace{.7cm}

\vskip1,2cm

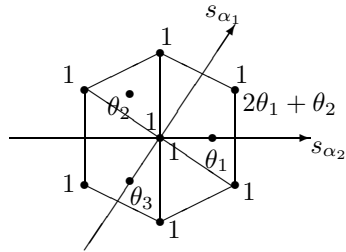
\begin{picture}(20,4)
\setlength{\unitlength}{1cm}
\put(4.7,0.73){$1$}
\put(5,0.63){\circle*{0.1}}
\put(5,0.63){\line(2,1){1}}
\put(5,0.63){\line(0,-1){1.26}}
\put(6,0){\line(3,-2){1}}
\put(6,0){\line(-3,2){1}}
\put(4.7,-0.73){$1$}
\put(5,-0.63){\circle*{0.1}}
\put(5,-0.63){\line(2,-1){1}}
\put(5.6,-0.58){\circle*{0.1}}
\put(5.6,-0.9){$\theta_3$}
\put(5.6,0.58){\circle*{0.1}}
\put(5.3,0.3){$\theta_2$}
\put(6,0){\circle*{0.1}}
\put(6,1.23){$1$}
\put(6,1.13){\circle*{0.1}}
\put(6,1.13){\line(0,-1){2.26}}
\put(6,0){\line(-1,0){2}}
\put(6,0){\vector(1,0){2}}
\put(6,0){\vector(2,3){1}}
\put(6,0){\line(-2,-3){1}}
\put(5.8,0.1){$1$}
\put(6.1,-0.3){$1$}
\put(6,-1.13){\circle*{0.1}}
\put(6.1,-1.33){$1$}
\put(6.7,0){\circle*{0.1}}
\put(6.6,-0.4){$\theta_1$}
\put(6.6,1.6){$s_{\alpha_1}$}
\put(7,0.73){$1$}
\put(7.1,0.35){$2\theta_1+\theta_2$}
\put(7,0.63){\circle*{0.1}}
\put(7,0.63){\line(-2,1){1}}
\put(7,0.63){\line(0,-1){1.26}}
\put(7.1,-0.83){$1$}
\put(7,-0.63){\circle*{0.1}}
\put(7,-0.63){\line(-2,-1){1}}
\put(8,-0.2){$s_{\alpha_2}$}
\end{picture}

\vskip1cm
\vspace{.7cm}
\noindent
For instance for $V^1 \otimes V^2$, the trivial small orthocells define in the
dual of the cartan subalgebras the 6 vectors, image of the weight $\varpi_1+\varpi_2$ under the action of the Weyl group (the vertices of the hexagon) and the 2 nontrivial correspond to two representations of a subgroup $SL(2,\BbbC)$ inside $SL(3,\BbbC)$ (two of the diagonals of the hexagon), thus to 2 times the weight 0, the third diagonal corresponding to a non small orthocell is excluded.\\
Let us now choose invariant supplementary space for $V^{ij}$ in
$V^i\otimes V^j$:
$$
\begin{array}{c}
V^1\otimes V^1 =V^{11}\oplus \vec \{
e_1\otimes e_2-e_2\otimes e_1; e_2\otimes e_3-
e_3\otimes e_2;\hfill\vspace{.2cm}\\
\hfill e_1\otimes e_3-e_3\otimes e_1\}\vspace{.2cm}\\
V^2\otimes V^2 =V^{22}\oplus \vec\{
(e_1\wedge e_3)\otimes (e_2\wedge e_3)-
(e_2\wedge e_3)\otimes (e_1\wedge e_3);\hfill\vspace{.2cm}\\
\hfill (e_1\wedge e_2)\otimes (e_1\wedge e_3)-
(e_1\wedge e_3)\otimes (e_1\wedge e_2);\vspace{.2cm}\\
\hfill (e_1\wedge e_2)\otimes (e_2\wedge e_3)-
(e_2\wedge e_3)\otimes (e_1\wedge e_2)\}\vspace{.2cm}\\
V^1\otimes V^2 =V^{12}\oplus \vec\{e_1\otimes (e_2\wedge e_3)+
e_2\otimes (e_3\wedge e_1)+e_3\otimes (e_1\wedge e_2)\}\vspace{.2cm}\\
V^2\otimes V^1 =V^{21}\oplus \vec\{(e_2\wedge e_3)\otimes e_1+
(e_3\wedge e_1)\otimes e_2+(e_1\wedge e_2)\otimes e_3\}\vspace{.2cm}\\
\end{array}
$$
We refind the sixteen relations defining the classical shape algebra
 by considering the relations:
$$
(I)_{ij}={\ker m}_{\vert{}_{\displaystyle V_i\otimes V_j}}
$$
$$
(II)_{12} = R_{12}(e_C^{12})=e_C^{21}, ~~ (II)_{11} ~~{\rm and}~~ (II)_{22}~~ {\rm are~~ trivials.}
$$
where $R_{12}$ is the isomorphism
$R_{12}: V^1\otimes V^2\rightarrow V^2\otimes V^1$.
Explicitly, we find :
$$
\begin{array}{cc}
(I)_{11}=\{ p_ip_j =p_jp_i, \hskip0,3mm (i<j)\}\hskip0,7cm\hfill
&\hfill(I)_{22}=\{ q_iq_j =q_jq_i, \hskip0,3mm (i<j)\}\vspace{.2cm}\\
(I)_{12}=\{ p_1q_1 +p_2q_2+ p_3q_3=0\}\hskip0,4cm
&(II)_{12}=\{ p_iq_j =q_jp_i, \hskip0,3mm \forall i,j\}\\
\end{array}
$$
To define the quantum $SL(3, \BbbC)$, we start with its representation theory,
 similar to the representation theory for classical $SL(3, \BbbC)$.
For instance, $V^1$ and $V^2$ becomes \cite{r2}:
$$
\begin{array}{c}
K_\beta e^i_{w}=q^{<w\varpi_i,\beta>}e^i_{w}\hfill\\
K_\beta^{-1} e^i_{w}=q^{-<w\varpi_i,\beta>}e^i_{w}\hfill\\
\end{array}
$$
and
$$
\begin{array}{ccc}
X_\beta e^i_{w}=0,\hskip3mm \hfill&Y_\beta e^i_{w}=e^i_{s_\beta w}
&\mbox{ if }<w\varpi_i,\beta>=1\hfill\vspace{.2cm}\\
X_\beta e^i_{w}=0,\hskip3mm\hfill & Y_\beta e^i_{w}=0\hfill
&\mbox{ if }<w\varpi_i,\beta>=0\hfill\vspace{.2cm}\\
 X_\beta e^i_{w}=e^i_{s_\beta w},\hskip3mm& Y_\beta e^i_{w}=0\hfill
&\mbox{ if }<w\varpi_i,\beta>=-1\\
\end{array}
$$
 for $(i,j)=(1,2)$ and $\beta =\alpha_1~~{\rm or}~~ \alpha_2.$\\
Now we define the quantum $e^{ij}_C$ by:
$$
e^{ij}_C =\sum_{L\subset A} q^{\vert L\vert}e^i_{s_{\bar L}w}
\otimes e^j_{s_{ L}w}.
$$
$V^{ij}$ is still generated by the $e^{ij}_C$, irreducible,
 with the highest weight $\varpi_i+\varpi_j$. We choose supplementary spaces
 for $V^{ij}$ in $V^i\otimes V^j$ by :
$$
\begin{array}{c}
V^1\otimes V^1 =V^{11}\oplus \vec \{(p_ip_j -qp_jp_i)^*,
 \hskip0,3mm (i<j)\}\vspace{.2cm}\\
V^2\otimes V^2 =V^{22}\oplus \vec\{(q_iq_j -qq_jq_i)^*,
 \hskip0,3mm (i<j)\}\vspace{.2cm}\\
V^1\otimes V^2 =V^{12}\oplus
 \vec\{(q^{-2}p_1q_1 +q^{-1}p_2q_2+ p_3q_3)^*\}\vspace{.2cm}\\
V^2\otimes V^1 =V^{11}\oplus \vec\{(q_1p_1 +q^{-1}q_2p_2+ q^{-2}q_3p_3)^*\}\\
\end{array}
$$
Now the sixteen relations defining the quantum shape algebra are :
$$
\begin{array}{cc}
&(I)_{11}\hskip8mm p_ip_j -qp_jp_i=0, \hskip0,3mm (i<j)\hfill\vspace{.2cm}\\
&(I)_{22}\hskip8mm q_iq_j -qq_jq_i=0, \hskip0,3mm (i<j)\hfill\vspace{.2cm}\\
&(I)_{12}\hskip8mm p_2q_2+q^{-1}p_1q_1 +q p_3q_3=0\hfill\vspace{.2cm}\\
&(I)_{21}\hskip8mm q_2p_2+qq_1p_1 + q^{-1}q_3p_3=0\hfill\vspace{.2cm}\\
&(II)_{12}=(II)_{21}\hskip8mm p_ip_j -qp_jp_i=0 \hskip2mm i\not= j\hfill\vspace{.2cm}\\
&p_1p_1 =qq_1p_1,\hskip2mm p_3q_3=q^{-1}q_3p_3.
\end{array}
$$
Since the quantum shape algebra is quadratic, we get here a description of this algebra.
\begin{theo}[Quantum shape algebra for ${\bf SL(3,\BbbC)}$]
\

\noindent
The quantum shape  algebra is the quotient of the tensor associative algebra
 $T(p,q)$ generated by $p_1,p_2,p_3,q_1,q_2,q_3$
 by the two sided ideal generated by these sixteen relations.
\end{theo}
\begin{rema}
In {\rm \cite{r1}} C. Ohn describe geometrically the preceding construction as a
 deformation of a sheme $E$  canonically defined in $\gb$.
\end{rema}

\section{Geometrical construction for $G_1$ and $G_0$}
In this part, we try to adapt the preceeding construction for the cases of
 $G_0$ and $G_1$.
We denote by $B_0$ and $B_1$ the ``Borel subgroup'' $B\cap G_0, ~B\cap G_1$
 ~for $G_0$ and $G_1$. We associate to them the  ''flag'' manifolds
$D_0=\gbo$ and $D_1=\gbun$. As in section \ref{sec} we get :
$$
\begin{array}{c}
D_0\simeq\{[p],[q]\in \BbbP(\BbbC^3),~~ q_3=1, ~~pq=0\}\subset D\vspace{.2cm}\\
D_1\simeq\{[p],[q]\in \BbbP(\BbbC^3), ~~q_3\not=0, ~~pq=0\}\subset D.
\end{array}
$$
Then $D_0=D_1$ is a dense subset of $D$. The characters for $D_0$ and $D_1$
 have the form :
$$
\begin{array}{c}
\chi_n^0:\left(\begin{array}{ccc}
a_1&b_1&0\\
0&b_2&0\\
0&0&1\\
\end{array}\right)
\mapsto a^{-n}_1\hskip3mm ~~(n\in \BbbZ)\vspace{.2cm}\\
\chi_{n_1,n_2}^1:\left(\begin{array}{ccc}
a_1&b_1&0\\
0&b_2&0\\
0&0&c_3\\
\end{array}\right)
\mapsto a^{-n_1}_1 c^{n_2}_3 \hskip3mm ~~(n_1,n_2\in \BbbZ).\\
\end{array}
$$
Thus ~~${\chi_{n_1,n_2}^1}\vert_{\displaystyle B_0}=\chi_{n_1}^0$~~ and the line
 bundles ~~$G_1\times_{\chi_{n_1,n_2}^1}\BbbC\rightarrow D_1$ ~~and
~~$G_0\times_{\chi_{n_1}^0}\BbbC\rightarrow D_0$ are isomorphic for any
$n_2$ \cite{r3}.\\
Let us consider now  the space of holomorphic sections for these  bundles.
 We put $U_i=G_i\cap U$ then :
$$
\begin{array}{c}
\guun\simeq\{(p,q)\in \BbbC^6,~~p\not=0 ~~q_3\not =0, ~~pq=0\}\vspace{.2cm}\\
\guo\simeq\{(p,q)\in \BbbC^6,~~p\not=0 ~~q_3=1, ~~pq=0\}.
\end{array}
$$
As in section \ref{sec} we embed these spaces in $\BbbP(\BbbC^7)$ and take their closure :
$$
\overline{\guun}=\overline{\gu}=
\{
\left[
\begin{array}{c}
p\\
q\\
1\end{array}
\right]
,\hskip3mm  p.q=0\}
$$
$$
\overline{\guo}=
\{
\left[
\begin{array}{c}
p\\
q\\
1\end{array}
\right]
,\hskip3mm q_3=1,  p.q=0\}\subset\overline{\guun}
$$
but ~~~~$\overline{\guo} \neq \overline{\guun}$
\begin{theo}[Space of section]

\

\noindent
$\bullet$
If $n_1<0$,
$$
H^0(\gbun,\chi_{n_1,n_2}^1)\simeq
H^0(\gbo,\chi_{n_1}^0)=~\{0\}.
$$
$\bullet$
If $n_1\geq 0$, then $H^0(\gbun,\chi_{n_1,n_2}^1)$ is infinitely dimensional.
\par
More precisely :
$$
H^0(\gbun,\chi_{n_1,n_2}^1)\simeq
{\displaystyle \large\bigcup_{l=sup(0,-n_2)}^\infty} \frac{1}{q_3^l}H^0(\gb,\chi_{n_1,l+n_2})
$$
and
$$
{\displaystyle \large\oplus_{n_1=0}^\infty}
{\displaystyle \large\oplus_{n_2=0}^\infty} H^0(\gbun,\chi_{n_1,n_2}^1)=
{\cal O}(\overline{\guun}\cap \{q_3\neq0\})
$$
Similarly :
$$
H^0(\gbo,\chi_{n_1}^0)\simeq
{\displaystyle \large\bigcup_{l=0}^\infty} H^0(\gb,\chi_{n_1,l})\vert_{q_3=1}
$$
and
$$
{\displaystyle \large\oplus_{n_1=0}^\infty}
 H^0(\gbo,\chi_{n_1}^0)\simeq
{\cal O}(\overline{\guo})
$$

\end{theo}
\vspace{.3cm}
{\bf Proof :} First, the space of holomorphic sections vanishes if $n=n_1<0$
since their restriction to $SL(2,\BbbC)$ are sections of the usual line bundle :
$$
SL(2,\BbbC)\times _{\chi_{n_1,n_2}}\BbbC\longrightarrow
{}^{\hbox{$SL(2,\BbbC)$}}\big/_{\hbox{$B\cap SL(2,\BbbC)$}}
$$
There are no restrictions on $n_2$.\par
As above, a section $f$ in $H^0(\gbun,\chi_{n_1,n_2}^1) $
can be viewed as a homogeneous function in
$x,y,z$ with degree $n_1+n_2$ in $x$, $n_2$ in $y$ and $0$ in $z_3$.
$f\in H^0(\gbun,\chi_{n_1,n_2}^1)$ has the form ~~
$ \displaystyle \frac{\varphi(x,y,z)}{(\det g)^k(q_3(g))^l}$. \\
We multiply $f$ by $(\det g)^k$ and we choose $l$ as  small as possible.\\
We get :
$$
 f = \displaystyle \frac{\varphi(x,y,z)}{(q_3(x,y))^l}.
$$
Moreover, the covariance relation $f(g_1b_1)=\chi_{n_1,n_2}^1(b_1)^{-1}f(g_1)$
implies that $\varphi$ is an homogeneous polynomial function with degree
 $0$ in $z$, $l$ in $y$ and $n_1+n_2$ in $x$.
$\varphi$ being $U_1$-invariant, since :
$$
$$$$
\left(\begin{array}{ccc}
x_1&y_1&0\\
x_2&y_2&0\\
x_3&y_3&z_3\\
\end{array}\right)
=\begin{array}{l}
\left(\begin{array}{ccc}
p_1&0&0\\
p_2&\frac{q_3}{p_1}&0\\
p_3&\frac{-q_2}{p_1}&\frac{1}{q_3}\\
\end{array}\right).
\left(\begin{array}{ccc}
1&\frac{y_1}{x_1}&0\\
0&1&0\\
0&0&1\\
\end{array}\right)
\end{array}.
$$
If $x_1\not= 0$, we can write :
$$
\begin{array}{cc}
f(x,y)&= \displaystyle \frac{F_1(p,q)}{p_1^{n_2}q_3^l}
\hskip9mm\hbox{for} ~~p_1\not = 0\vspace{.2cm}\\
\hfill&= \displaystyle \frac{F_2(p,q)}{p_2^{n_2}q_3^l}
\hskip9mm\hbox{for} ~~p_2\not = 0\\
\end{array}
$$
As in the section \ref{sec}, the function ~$q_3^l f$ coincides in fact with a polynomial function
 in the variables $p$ and $q$. Thus :
$$
(**)~~~f(x,y)= \displaystyle \frac{1}{q_3^l}F(p,q)
\hskip3mm(\deg_p F =n_1, \deg_qF= l+n_2\geq 0)
.$$
But now ~~$ 1/{q_3^l}$~~is in ${\cal O}(\gu)$, we can not eliminate
the denominator $q_3^l$.\\ By the preceeding
 discussion, $F$ can be viewed as an  element of
$H^0(\gb,\chi_{n_1,l})$ and $(**)$ proves that :
$$
H^0(\gbun,\chi_{n_1,n_2}^1) ~\simeq
 ~{\displaystyle \large\bigcup_{l=sup(0,-n_2)}^\infty} \frac{1}{q_3^l}H^0(\gb,\chi_{n_1,l+n_2})
$$
We don't have a direct sum since for instance:
$$
1 = \frac{q_3}{q_3} = ... =  \frac{q_3^l}{q_3^l} \in H^0(\gb,\chi_{0,0})\large\cap \frac{1}{q_3}H^0(\gb,\chi_{0,1})\large\cap ...\large\cap \frac{1}{q_3^l}H^0(\gb,\chi_{0,l})
$$
Finally, $\overline{\guun} = \overline{\gu}$, then a function $f$ belongs to ${\cal O}(\overline{\guun}\cap \{q_3\neq0\})$ if and only if $f$ can be written:
$$
f(x,y)= \displaystyle \frac{1}{q_3^l}F(p,q)~~~~~~~(l\geq0),$$
where $F$ is a polynomial function. Then $F = \sum_{n_1,n'_2}F_{n_1,n'_2}$ where $F_{n_1,n'_2}$is a homogeneous with degree $n_1$ in $p$, $n'_2$ in $q$ or:
$$
f(x,y)= \displaystyle \sum_{n_1,n'_2\geq 0}\frac{1}{q_3^l}F_{n_1,n'_2}= \sum_{n_1\geq 0,n_2\geq -l}\frac{1}{q_3^l}F_{n_1,n_2+l} \in \sum_{n_1,n_2}H^0(\gbun,\chi_{n_1,n_2}^1)
.$$
The sum is clearly direct.
Similarly, since the bundle $G_{0}\times _{\chi_{n_1}^0}\BbbC$ is the restriction to $G_0$ of the bundle $G_{1}\times _{\chi_{n_1,0}^1}\BbbC$ (\cite{r3}), we have:
$$
H^0(\gbo,\chi_{n_1}^0)= H^0(\gbun,\chi_{n_1,0}^1)\vert_{q_{3}=1}
$$
this gives the last assertions of our theorem.

\section{Shape algebra for $G_1$}

Similarly to the $G$-case, we shall define the (classical and quantum) shape algebras for $G_1$ as the vector space ${\cal O}(\overline{\guun}\cap \{q_3\neq0\})$ of all the line bundles over $\gbun$.\\
At the classical level this algebra is generated by the space $V_1\oplus V_2\oplus V_{-1}$ where:
$$
V_1 = Vec(p_1, p_2, p_3),~~~~V_2 = Vec(q_1, q_2, q_3), ~~~~V_{-1} = \BbbC \frac{1}{q_3}
$$
In order to define the multipication law of our shape algebra, we need a large family of representations $(\pi_i)$ of $G_1$ such that the representation $\pi_i\otimes\pi_j$ is a finite sum of some $\pi_k$.\\
Let us recall that $G_1$ is a classical and quantum subgroup of $G$ (see section 2), then each irreducible finite dimensional representation $V^\lambda$ of $G$ is a representation, still denoted $V^\lambda$ of $G_1$.\\
Since $V^\lambda$ is generated by a highest weight vector (and its dual $V_\lambda$ by a lowest weight vector, the funtion $p_3^{n_1}q_1^{n_2}$), then $V^\lambda$ is an indecomposable representation of $G_1$, generally it is not irreducible: the $SL(2)$ module, generated by the highest weight vector, is a $G_1$-submodule without any direct factor.\\
We select now, the family ~$((V^{-1})^{\otimes l}\otimes V^{\lambda})_{n_1\geq0, n_2\geq0, l\geq0}$ ~as our family of representation for $G_1$ the dual of such a representation appears naturally as the space of functions $f$ in ${\cal O}(\overline{\guun}\cap \{q_3\neq0\})$ of the form:
$$
f(x,y)= \displaystyle \frac{1}{q_3^l}F(p,q)~~~~deg_pF = n_1, ~~deg_qF = n_2.
$$
Especially it is generated by $V_{-1}$, $V_1$ and $V_2$, since $V_{-1}\otimes V_{i}$ and $V_{i}\otimes V_{-1}$ are elements of our family of representation, there is no supplementary spaces, thus no relation like $I_{-12}$ or $I_{-11}$. However, the shape algebra is no more a direct sum of our representation spaces.Thus we have to add a new relation $I^0_{-1,2} = I^0_{2,-1}$ since
$$
(V_{-1})^{\otimes 0} \otimes V_{0,0} \subset V_{-1}\otimes V_{2}$$
$$
1 \otimes 1 = \frac {1}{q_3} \otimes q_3 = q_3 \otimes \frac {1}{q_3}
$$
Thus the quadratic relations for the classical shape algebra of $G_1$ are those of the classical $G$:\\
$$
I_{11}, ~~I_{12} = I_{21},~~ I_{22}, ~~II_{12} = II_{21} $$
and
$$
\begin{array}{c}
I^0_{-12}~~~~ (1/{q_3}) . q_3 = 1\\
I^0_{2-1}~~~~ q_3 . (1/{q_3}) = 1\\
II_{-11}~~~~(1/{q_3}) . p_i = p_i . (1/{q_3})\\
II_{-12}~~~~(1/{q_3}) . q_i = q_i . (1/{q_3})
\end{array}
$$
But the relations $II_{-11}$ and $II_{-12}$ are consequences of $I^0_{-12}$, $I^0_{2-1}$ and
$I_{ij}$, $II_{ij}$  $(i\geq0, j\geq0)$.\\
Let us, now consider the quantum case, thus $V^1$, $V^2$ are the restriction to $G_1$ of the $G$ module $V^1$, $V^2$, we define $V^{-1}$ as the one dimensional space $\BbbC v$ with:
$$
X_1 v = Y_1 v = Y_2 v = 0, ~~~~~~~~K_1 v = v, ~~~~~~K_2 v = q^{-1}v
.$$
Then the quadratic relations in the shape algebra are:
$$
I_{11}, ~~I_{22},~~I_{12},~~ I_{21}, ~~II_{12} = II_{21} $$
and
$$
\begin{array}{c}
I^0_{-12}~~~~ (1/{q_3}) . q_3 = 1\\
I^0_{2-1}~~~~ q_3 . (1/{q_3}) = q\\
\end{array}
$$
\begin{theo}[Quantum shape algebra for ${\bf G_1}$]

\

\noindent
The quantum shape algebra is the quotient of the tensor associative algebra $T(p,q,1/{q_3})$, generated by $p_1, p_2, p_3, q_1, q_2, q_3, 1/{q_3}$, by the above relations.
\end{theo}
Indeed, the shape algebra is a quotient of the algebra defined in theorem 3 but if $q = 1$ this quotient has to coincide with ${\cal O}(\overline{\guun}\cap \{q_3\neq0\})$ which is the classical shape algebra, thus our quotient is trivial.

\end{document}